\documentclass[a4paper, 12pt]{amsart}
\usepackage{graphicx}
\newtheorem{thm}{Theorem}[section]

\newtheorem{lem}[thm]{Lemma}

\theoremstyle{definition}
\newtheorem{defn}[thm]{Definition}
\theoremstyle{remark}

\numberwithin{equation}{section}

\newcommand{\set}[1]{\left\{#1\right\}}

\def\abra#1#2#3{\begin{figure}[ht]\begin{center}
\includegraphics[width=#3cm]{#1.eps}
 \end{center}
  \caption{#2}
  \label{#1}
\end{figure}}

\def\set#1{\{#1\}}

\begin{document}

\title[Graphs with path-width two]{On the structure of graphs with path-width at most two}

\author[J. Bar\'at]{J\'anos Bar\'at}
\address{Bolyai Institute, University of Szeged, Aradi v\'ertan\'uk tere 1, 6720 Szeged, Hungary}
\curraddr{Department of Computer Science, University of Pannonia, Egyetem u. 10, 8200 Veszpr\'em, Hungary}
\thanks{Research is supported by OTKA Grant PD~75837.}
\email{barat@dcs.vein.hu}

\author[P. Hajnal]{P\'eter Hajnal}
\address{Bolyai Institute, University of Szeged, Aradi v\'ertan\'uk tere 1, 6720 Szeged, Hungary}
\thanks{Research is supported by OTKA Grant K~76099.}
\email{hajnal@math.u-szeged.hu}

\author[Y. Lin]{Yixun Lin}
\address{Department of Mathematics, Zhengzhou University, Zhengzhou, Henan 450052, China}
\thanks{Supported by NSFC grant 10071076}
\email{linyixun@zzu.edu.cn}

\author[A. Yang]{Aifeng Yang}
\address{School of Management, Hefei University of Technology, Hefei 230009, China}



\dedicatory{Dedicated to Professor Carsten Thomassen
on the occasion of his sixtieth birthday}


\begin{abstract}
Nancy G.~Kinnersley and Michael A.~Langston
has determined the excluded minors for the class
of graphs with path-width at most two by computer.
Their list consisted of 110 graphs.
Such a long list is difficult to handle and gives no insight
to structural properties.
We take a different route, and concentrate on the building blocks
and how they are glued together.
In this way, we get a characterization of $2$-connected and
$2$-edge-connected graphs with path-width at most two.
Along similar lines, we sketch the complete characterization of 
graphs with path-width at most two.
\end{abstract}

\maketitle

\section{Introduction}

The concept of path-width and tree-width
were defined
and played central role in the Graph Minors
project by Robertson and Seymour. One important feature of path-width and
the first result of this type is that excluding a tree in an infinite
sequence of finite graphs results in a class of bounded path-width \cite{bie2, die2, gm1}.
Similarly, a forbidden planar graph implies that the class has bounded
tree-width \cite{gm4}.
These results are complemented with the theorem that graphs of bounded
tree-width are well-quasi-ordered \cite{gm5} giving a prototype of the deep and
lengthy Graph Minor Theorem.
The proof for bounded tree-width is comprehensible and now well digested.
This is one of the reasons that path-width is sometimes considered too
simple or less valuable.
Since bounded tree-width implies bounded path-width, there is no direct
proof that the graphs of bounded path-width are well-quasi-ordered.
If path-width is "so simple", then there should be a
Nash-Williams type proof for this latter result \cite{nash}.
However, we believe that such a result is still unknown.

Another dogma we would like to attack is that graphs of large path-width
are more important than that of small path-width.
Our point is that raising the connectivity and path-width simultaneously gives
structural information.
This idea has also been exploited in \cite{gup} for algorithmic use.
The number of excluded minors of this kind seems to increase mildly opposite to
the number of excluded trees, which grows super-exponentially \cite{kin,tak}.

\section{Notations}

We consider finite, simple graphs except in Section~\ref{2elof}, where allowing
double edges makes the list of excluded minors more compact.
A graph $H$ is a \sl minor \rm of a graph $G$, denoted as $G\succeq H$,  
if $H$ can be obtained from a subgraph of $G$ by contracting edges.
Contraction of an edge might lead to double edges.
When we consider 2-edge-connected graphs, it is natural to allow
double edges.
It will make our discussion more comfortable.
Otherwise, we keep the graph simple by removing multiple edges after contraction.

In this paper, we focus on the following well-known graph parameter.

\begin{defn}
A {\rm path-decomposition} of a graph $G$ is a pair $(P,W)$, where $P$
is a path, and $W=(W_p : p \in V(P))$ is a family of subsets of
$V(G)$, satisfying\\
\indent $(1)$
$\displaystyle
\bigcup_{p\in V(P)} W_p= V(G)$, and every edge of $G$ has both ends
in some $W_p$, and\\
\indent $(2)$
if $p,p',p''\in V(P)$ and $p'$ lies on the path from $p$ to $p''$,
then $\displaystyle W_p\cap W_{p''}\subseteq W_{p'}$.

The {\rm width}
of a path-decomposition is
$\max (|W_p|-1 : p\in V(P))$, and the
{\rm path-width}
of $G$ is the minimum width of all
path-decompositions of $G$.
\end{defn}

\abra{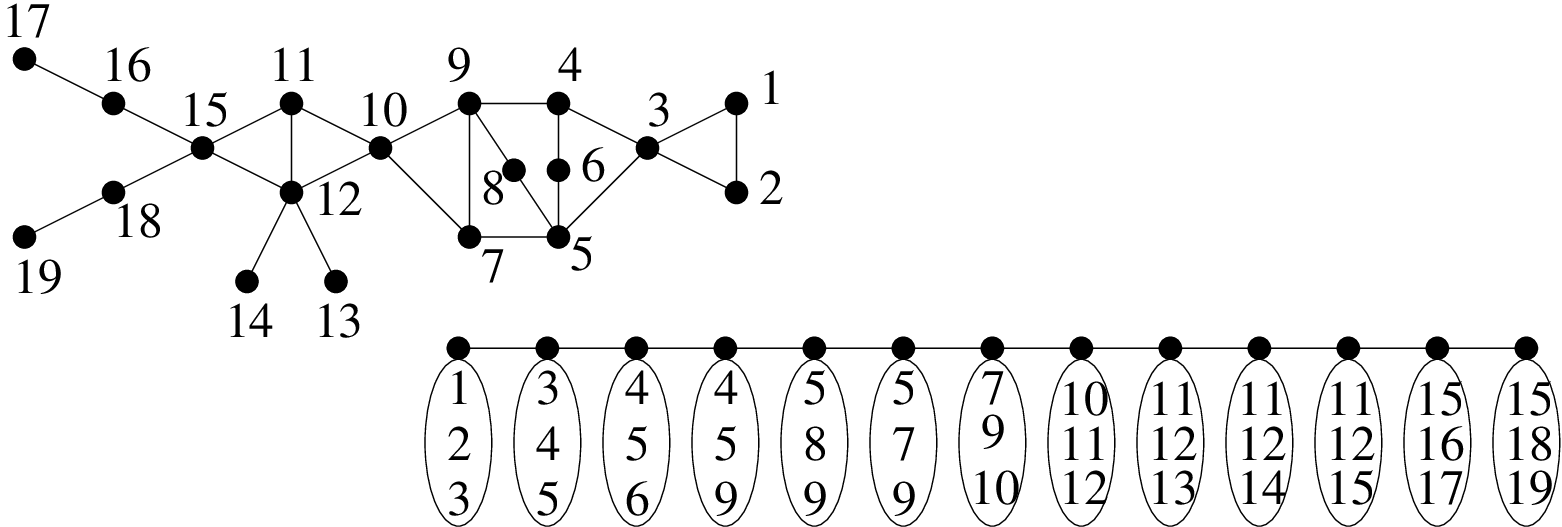}{Path-width two corresponds to bags of size three}{12}

It follows from the definition that path-width is minor-monotone,
that is, if $G\succeq H$, then $pw(G)\ge pw(H)$.
We will denote the class of graphs with path-width at most 2 by PW2.
It is clear that if $G\in$PW2, then $G\not\succeq K_4$.
Therefore, any graph in PW2 is a planar graph.

We will not make any use of it, but path-width is equivalent
to a cops-and-robber game \cite{tho}. Being familiar
with the game might help the reader's intuition throughout
the discussion.

A {\sl rooted graph} $(G,r)$ is a graph $G$ with a specific
node $r\in V(G)$, that is called the {\sl root} of $G$.
The rooted graph $(H,s)$ is a rooted minor of $(G,r)$,
denoted as $(G,r)\succeq_r (H,s)$, if $r$ is mapped to $s$,
when $G$ is mapped to $H$ under the minor operation.
Similarly, a two-rooted graph $(G,r_1,r_2)$ is a graph with
two specified nodes. Rooted minor of a two-rooted graph is
defined analogously. With a slight abuse of notation, $\succeq_r$
will be used for both rooted and two-rooted minors.

Let $G$ be a graph and let $U$ be a set of vertices.
Then $E(U)$ is the set of edges incident to any element of $U$.
The graph $G-U$ is obtained by deleting the vertices in $U$.
The graph $G|_U$ is the subgraph of $G$ induced by $U$.
Let $F$ be a set of edges.
Then $V(F)$ is the set of vertices
incident to at least one member of $F$.
The graph $G-F$ is obtained from $G$ by deleting
the edges in $F$.
The graph $G|_F$ is the subgraph induced by the edge set $F$,
i.e.~its vertex set is $V(F)$ and its edge set is $F$.

It is well-known, that ``being
identical or being on the same cycle''
is an equivalence relation
on the edges of a graph.
Its equivalence classes span
the so-called {\sl blocks} of the graph.
The one-element classes are the cut edges.

Let $G$ be a connected graph, and let $C$ be a cycle of $G$.
Consider $G-C$, and let $N$ be the vertices of a component of $G-C$.
Then $G|_{E(N)}$ is a {\sl bridge} of $C$ (in $G$), as well as any
edge connecting two nodes of $C$.
The {\sl legs of a bridge $B$} are the common vertices of $B$ and $C$.
The set of legs are denoted by $L(B)$.
Let $\set{a_1,a_2}$ and $\set{b_1,b_2}$
be two-element subsets of $V(C)$.
The two pairs are {\sl crossing} if and only if they are disjoint, and
along the cycle, the $a$-vertices alternate with the $b$-vertices.
Let $U$ and $V$ be subsets of $V(C)$.
We say that $U$ and $V$ are {\sl crossing}, if there are two nodes of $U$
and two nodes of $V$ such that the two pairs are crossing.
{\sl Two bridges are crossing} if and only if their set of legs are crossing.
A bridge is {\sl simple} if and only if it is a path (and  hence it
has two legs, the end-vertices of the path).
A bridge is called trivial, if it consists of one edge or two edges.

We describe a special class of graphs called tracks. They turn out to have
path-width at most two and play fundamental role in our discussion.

Let $P=\set{p_1,p_2,\ldots,p_k}$ and $Q=\set{q_1,q_2,\ldots,q_\ell}$ be two vertex disjoint paths.
This is a slight abuse of notation, sometimes $k=1$ and $P$ is a path of length $0$.
We define two types of connections between the two paths.
Firstly, there might be edges connecting a vertex of $P$ to a vertex of $Q$.
We call these edges {\sl short chords}.
Secondly, we allow paths of length two connecting $P$ and $Q$.
We call these paths {\sl long chords}.
A long chord has three nodes, a $p_i$, a middle node $m$, and
a $q_j$. The degree of $m$ is $2$ in the graph $G$.
That is, different long chords must have different
middle nodes, in particular long chords are edge-disjoint.
If we do not want to specify whether we talk about a short- or a long chord we refer
to it as a chord.
For brevity, we call a $p_iq_j$-chord an $ij$-chord.
An $ij$- and an $i^\prime j^\prime$-chord are not crossing if and only if $(i-i')(j-j')\geq 0$.
A $11$-chord or a $k\ell$-chord is called an {\sl end-chord}.
If there are several candidates for an end-chord, then one of them has to be selected.

\begin{defn}
A graph $G$ is called a
{\sl track} if and only if it can be represented
by two vertex disjoint paths $P, Q$ and
noncrossing chords as above such that $p_1$ and $q_1$ are connected by a chord
as well as $p_k$ and $q_\ell$.
\end{defn}

In this way, a track is set to be $2$-connected.

Alternatively, we can look at a track as a graph that can
be obtained from special outerplanar graphs by certain operations.

Let $G$ be an outerplanar graph with outer cycle $C$.
A graph $G$ is called {\sl multichordal} outerplanar if we allow multiple
edges inside $C$.
If we have a drawing of such a graph, then there might be some $2$-faces.
If the dual of the interior of $G$ is a path, then $G$ is called {\sl series}.
Consider the simple graphs arising from a series multichordal outerplanar 
graph by subdividing some of its internal edges by a single vertex.

These graphs are tracks for the following reason.
We can select two paths $ab$ and $a'b'$ of the outer cycle $C$ such that

\begin{itemize}
 \item each of $aa'$ and $bb'$ is an edge or a path with two edges;
 \item the cycle $C=abb'a'a$;
 \item all chords have one end on the path $ab$ and the other end on $a'b'$.
\end{itemize}

\begin{figure}[ht]
	\begin{center}
	\includegraphics[scale=0.5]{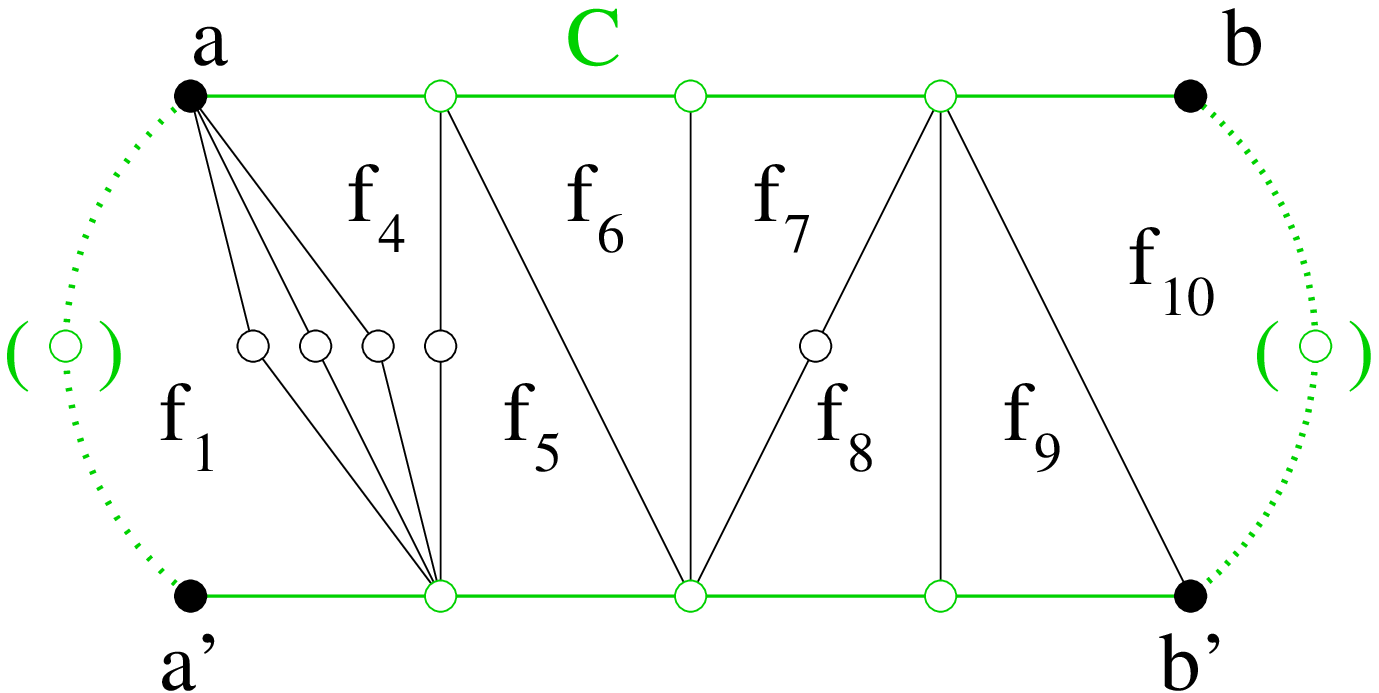}
 	\end{center}
  \caption{A track with its corners}
\end{figure}

The vertices $a,a'b,b'$ (that is $p_1, p_k, q_1, q_{\ell}$ in the original definition) 
will be called {\sl corners}.
We imagine $a$ and $a'$ to be on the left-hand side,
$b$ and $b'$ on the right-hand side.
Therefore, $a$ and $b$ will be called {\sl opposite} corners,
as well as $a'$ and $b'$.
In the general case, a track has four corners.
But some of them may coincide.
If $a=b$, then the node $a$ is called a degenerate side of the track.

If we have a drawing of a track, then the inner faces can be listed as
$f_1,f_2,\dots f_k$ according to the dual path.
Then $a$ and $a'$ are on face $f_1$ and $b$ and $b'$ are on face $f_k$.  
Notice that the selection of the corners are not unique.
Also the drawing of a track is not unique. 
This can make our proofs more complicated, but only technically.

\begin{defn}
A graph $G$ is called a partial track if and only if it is a
subgraph of a track.
\end{defn}



\section{Two-connected graphs}

It follows from the definition that tracks are 2-connected and
have path-width at most two.
We put this observation in a broader context.

\begin{thm} \label{2conntrack}
The following three statements are equivalent.
\begin{enumerate}
\item[(i)]
The graph $G$ is $2$-connected and $pw(G)\le 2$.
\item[(ii)]
The graph $G$ is a track.
\item[(iii)]
The graph $G$ is $2$-connected and $G\not\succeq E_1, E_2, E_3$.
\end{enumerate}
\end{thm}

\begin{figure}[ht]
	\begin{center}
	\includegraphics[scale=0.5]{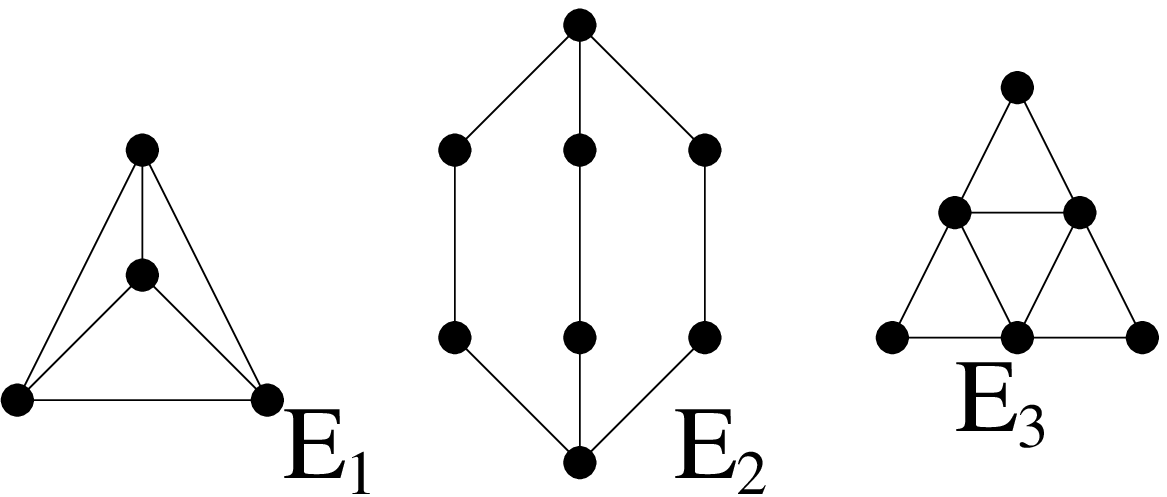}
	\end{center}
  \caption{Excluded minors for two-connected PW2-graphs}
\end{figure}

\begin{proof}
Two implications are immediate: $(i)\Rightarrow (iii)$ and $(ii)\Rightarrow (i)$.
We show that $(iii)$ implies $(ii)$.
In a 2-connected graph any two vertices are contained in a cycle.
Let $C$ be a longest cycle of $G$.
Consider the bridges of $C$.
Since $G\not\succeq E_1$, there is no $3$-bridge.
Since $G\not\succeq E_2$, there is no internal edge in any bridge.
Therefore the bridges are trivial.
Since $G\not\succeq E_1$, the bridges are pairwise equivalent or avoiding.
Since $G\not\succeq E_3$, the graph $G$ is a track.
\end{proof}

Note that the proof gives guidelines to find $a,b,a',b$ but usually
these four vertices are not well-defined. Even if the proof points
out certain vertices, the chords, that play the role of end-chords,
might be ambiguous.


\section{Two-edge-connected graphs} \label{2elof}

Our goal is to prove an excluded minor characterization of
2-edge-connected PW2-graphs.
We have to focus on the blocks, and how they are glued together.
Let $G\in$PW2 be 2-edge-connected.
We know from basic graph theory that $G$ is built up from its blocks pasted together
along vertices.
We call these attachment vertices {\sl multiple points}.

We collect some information about the blocks and the multiple points in four statements.
When we say that a node $r$ can be a corner of a track $T$, we literally mean that $r$ can
be fixed to be corner $a$, and the other three corners of $T$ can be selected in such
a way that we get back the definition of a track with corners $a,b,a',b'$.
Similarly, we say that $r$ can be a degenerate side of a track $T$ meaning that after fixing
$r=a=b$, we can select the other two corners properly.

The proofs are based on the following visualization of a track:
it is a 2-connected graph, consider its longest cycle $C$.
Since the chords are non-crossing bridges of $C$, they can be linearly
ordered.
This order is not well-defined, because equivalent bridges can be interchanged.
If the order is fixed, the position of $r$ is crucial,
and we focus on that.

\begin{lem} \label{s1}
Let $(T,r)$ be a track with a specified node.
If $T\not\succeq_r R_1, R_2$, then $r$ can be a corner of $T$.
\end{lem}

\begin{figure}[ht]
	\begin{center}
	\includegraphics[scale=0.5]{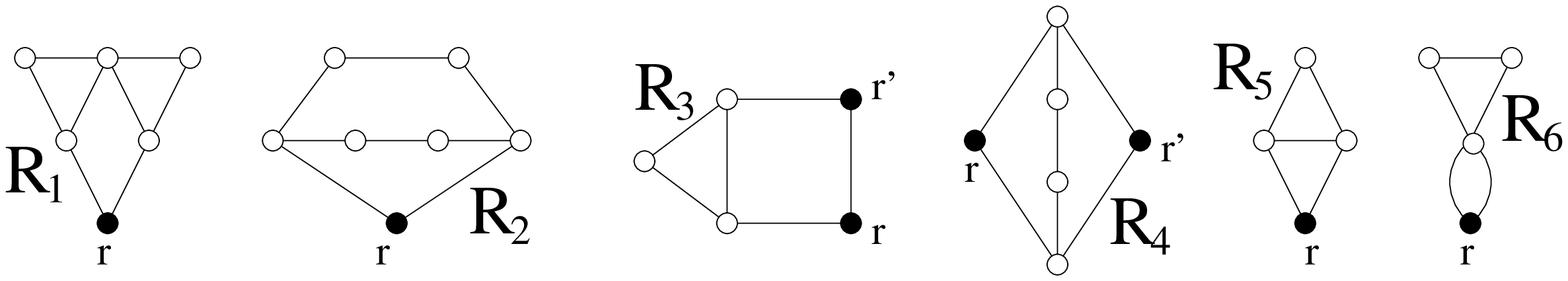}
 	\end{center}
  \caption{Crucial building bricks with specific gluing points}
\end{figure}

\begin{proof}
Let $C$ be a longest cycle through $r$ in $T$.
Consider the bridges of $C$.
Since $T$ is a track, there is no bridge with three legs.
Since $(T,r)\not\succeq_r R_2$, the bridges avoiding $r$ are chords.

Assume there is a bridge $H$ adjacent to $r$ and $s$.
If there was a path of length three from $r$ to $s$ in $H$, then
$T\succeq E_2$, a contradiction.
That means $H$ is a single vertex or empty.
The vertices $r$ and $s$ cuts the cycle $C$ into two parts,
let us say left and right from $r$.
If there is a chord avoiding $r$ on both sides, then $(T,r)\succeq_r R_1$.
Otherwise all chords are on one side, left from $r$ say.
We can now select the two paths $P$ and $Q$ in the definition of a track.
Let $r$ be a corner, $a$ say. The right neighbor of $r$ will be $a'$.
Starting a path from $r$ going to the left, we detect the legs of the chords in order.
There is a unique moment, when two consecutive legs $f$ and $g$ belong to the same chord.
If there were two such chords, then $T\succeq E_3$, a contradiction.
Therefore $f$ is a good choice for $b$ and the left neighbor of $f$ for $b'$.

The case when there are no bridges adjacent to $r$ is very similar.
Going on $C$ to the left and to the right from $r$, we find the legs $d$ and $e$
of the same chord. Continuing to the left from $d$ we detect another place
where two consecutive legs $f$ and $g$ belong to the same chord.
This is unique, since $T\not\succeq E_3$.
Now $r$ can play the role of $a$, the right neighbor of $r$ can be $a'$.
We select $f$ to play the role of $b$ and the left neighbor of $f$ can be $b'$.
\end{proof}

\begin{lem} \label{s2}
Let $(T,r,r')$ be a track with two specified nodes such that any of $r$ and $r'$ can be a corner of $T$.
If $T\not\succeq_r R_3, R_4$, then $r$ and $r'$ can be two opposite corners.
\end{lem}

\begin{proof}
Let $C$ be a longest cycle through $r$ and $r'$ in $T$.
Consider the bridges of $C$.
Since $T$ is a track, there is no bridge with three legs.
The vertices $r$ and $r'$ cut the cycle $C$ into two parts, let us say
left side and right side.
If there is a bridge of $C$ with two legs on the same side, then
$T\succeq_r R_3$, a contradiction.

If there is a bridge with legs $r$ and $r'$, then the bridge is trivial,
as otherwise $T\succeq E_2$.
Also all other bridges are trivial and adjacent to $r$ or $r'$.
Since $(T,r)\not\succeq_r R_1$, there are no two bridges adjacent to $r'$ such that
one of the other legs are on the left, and one is on the right from $r$.
Therefore, there might be bridges $H_1,\dots, H_k$ adjacent to $r$ with legs $s_1,\dots, s_k$
on the right side, and bridges $H'_1,\dots, H'_k$ adjacent to $r'$ with legs $s'_1,\dots, s'_k$ on the left side.
In this case $r$ and $r'$ can be corners as follows:
we select $r=a$ and $r'=b'$, and the path $P$ is the part of $C$ on the right side and
path $Q$ is the part of $C$ on the left side. Therefore, $a'$ is the left neighbor of $a$
and $b$ is the right neighbor of $b'$.

\begin{figure}[ht]
	\begin{center}
	\includegraphics[scale=0.5]{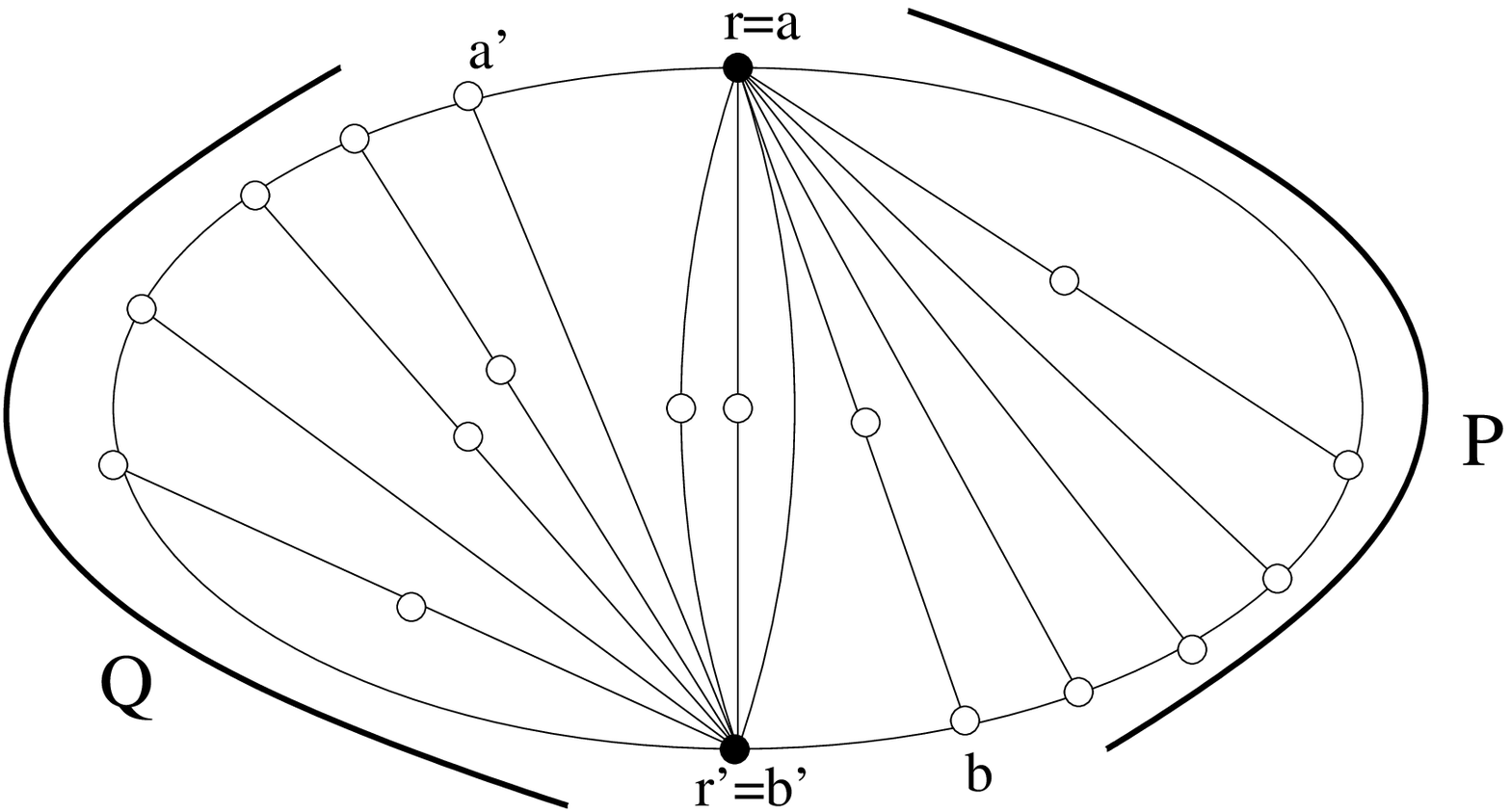}
 	\end{center}
  \caption{An illustration for the proof}
  \label{baljobb}
\end{figure}

If there is no bridge with legs $r$ and $r'$, then assume there is a bridge
with legs crossing $r$ and $r'$. Since $T\not\succeq_r R_4$, this bridge is trivial.
There might be many non-crossing chords like that.
Finally, there might be some chords on one side adjacent to $r$ and some on the
other side adjacent to $r'$. In this case the same scenario works as before.
We select $r=a$ and $P$ to be the left side of $C$, and $r'=b'$ and $Q$ to be
the right side of $C$.

\end{proof}

\begin{figure}[ht]
	\begin{center}
	\includegraphics[scale=0.5]{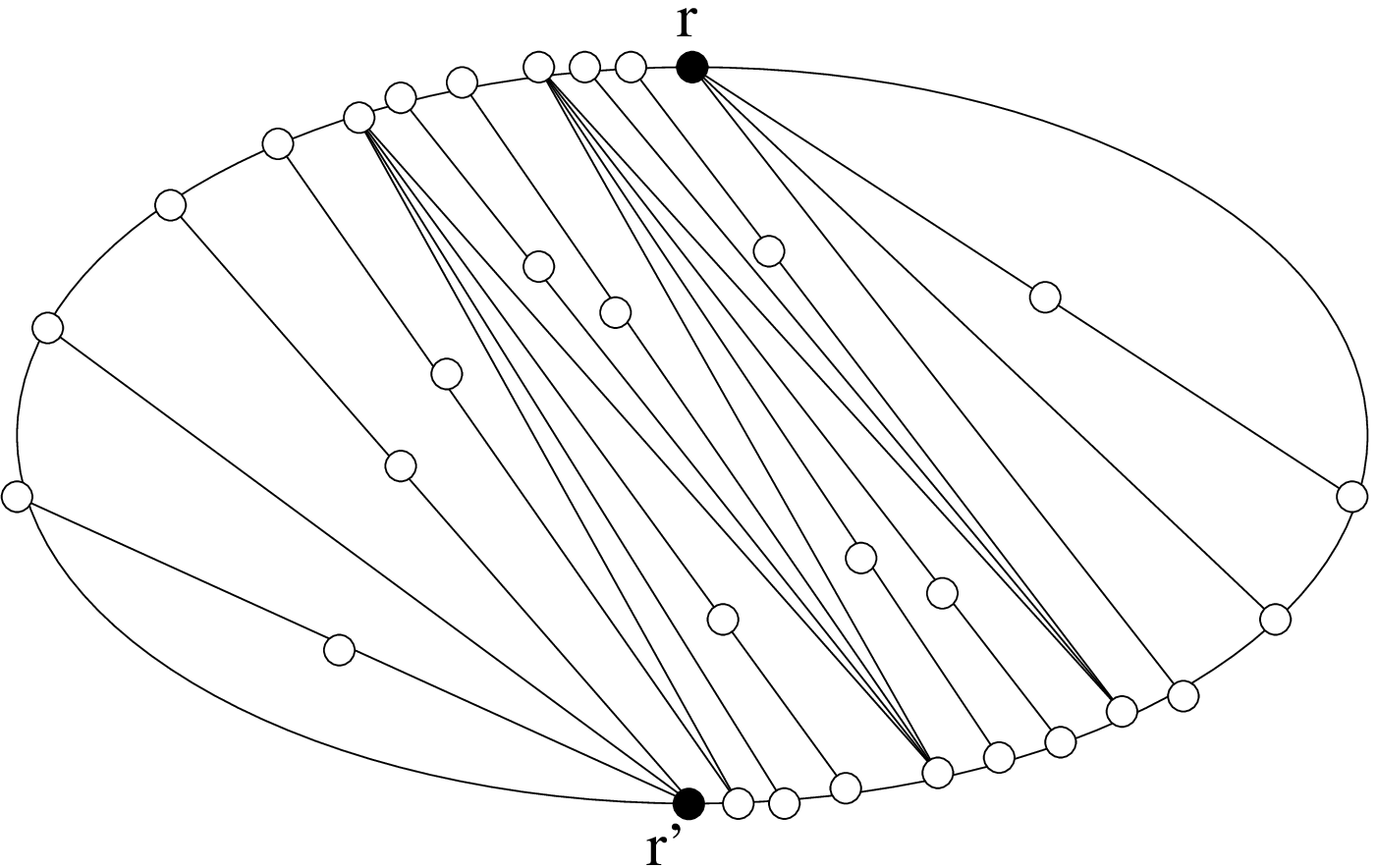}
 	\end{center}
  \caption{Another case of Lemma~\ref{s2}}
  \label{cross}
\end{figure}

\begin{lem} \label{s3}
Let $(T,r)$ be a track with a specified node.
If $T\not\succeq_r R_5$, then $r$ can be a degenerate side of $T$.
\end{lem}

\begin{proof}
Let $C$ be a longest cycle through $r$ in $T$.
Since $T\not\succeq_r R_5$, there is no bridge avoiding $r$.
If there is a nontrivial bridge with leg $r$ and other leg $s\in C$,
then the length of $C$ between $r$ and $s$ is at least three on both sides,
and therefore $T\succeq E_2$, a contradiction.
So $T$ consists of $C$ and some trivial bridges, all of them adjacent to $r$.
These graphs satisfy the claim.
\end{proof}

We also need the following observation.

\begin{lem} \label{s4}
Let $(T,r)$ be a $2$-edge-connected graph in which each block is a track, and 
$r$ is a specified node.
If $T\not\succeq_r R_6$, then the blocks of $T$ are pasted along $r$.
\end{lem}


We are done with all preparation for the main result, that is an excluded minor characterization
of $2$-edge-connected PW2-graphs.

\begin{figure}[ht]
	\begin{center}
	\includegraphics[scale=0.5]{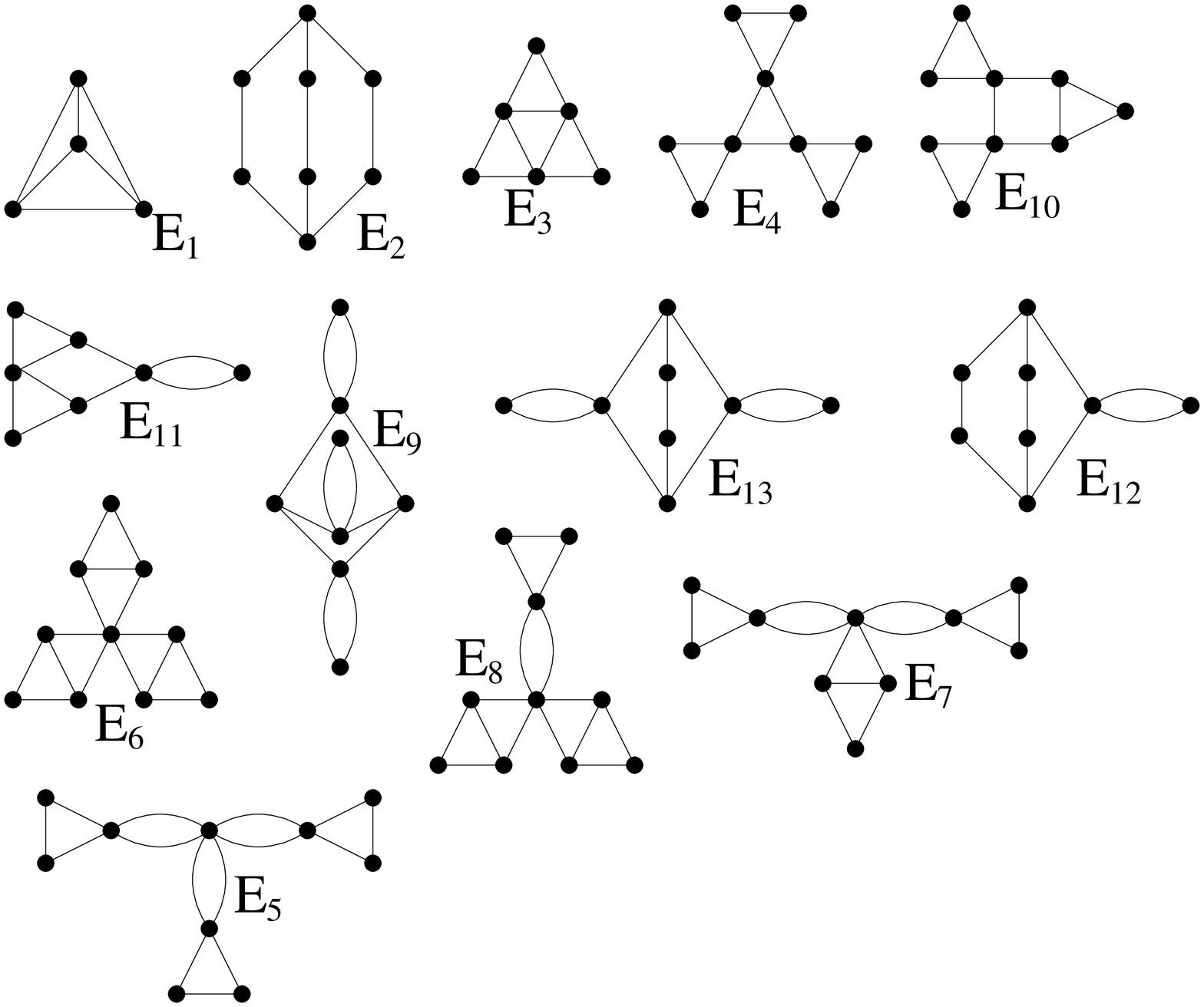}
 	\end{center}
  \caption{The excluded minors for $2$-edge-connected graphs}
  \label{excl}
\end{figure}

\begin{thm} \label{fo}
The following three statements are equivalent.
\begin{enumerate}
\item[(i)]
The graph $G$ is a $2$-edge-connected partial track.
\item[(ii)]
The $2$-edge-connected graph $G\not\succeq E_1, E_2,\dots, E_{13}$ shown in Figure~$\ref{excl}$.
\item[(iii)]
(The blocks of $G$ are tracks glued together in a path-like fashion according
to the following scheme)
The blocks of $G$ are tracks and can be listed as $B_1,B_2,\dots, B_L$
such that for any $2\le i\le L-1$ the block $B_i$ is a block with two (possibly identical) multiple nodes
$m_i$ and $m_{i+1}$.\\
If $m_i\neq m_{i+1}$, then $m_i$ and $m_{i+1}$ can play the role of opposite corners in $B_i$.\\
If $m_i=m_{i+1}$, then this multiple node can play the role of a degenerate side in $B_i$.\\
For $i=1$ or $L$ the block $B_i$ has a multiple node $m_i$ that can play the role of a corner.
\end{enumerate}
\end{thm}

\begin{figure}[ht]
	\begin{center}
	\includegraphics[scale=0.5]{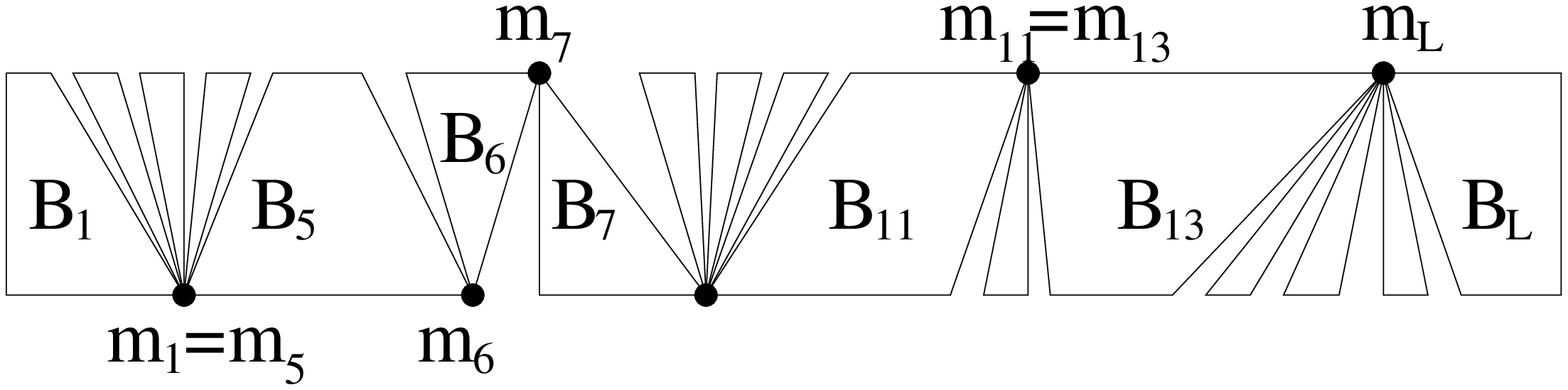}
 	\end{center}
  \caption{The gluing pattern of blocks}
  \label{glue}
\end{figure}

Notice that the above numbering of blocks is not unique.

\begin{proof}
Two implications are easy to see: $(i)\Rightarrow (ii)$ and $(iii)\Rightarrow (i)$.
We show that $(ii)$ implies $(iii)$.

In the block structure of $G$, the excluded minors $E_4$ and $E_9$ imply that 
no block can contain more than two multiple nodes.

The excluded minors $E_{11}$ and $E_{12}$ together with Lemma~\ref{s1} imply
that the block structure of $G$ satisfy that $m_i$ can play the role of a corner in $B_i$, for any $i$.

If there is a block $B$ with two distinct multiple points $m$ and $m'$, then
$E_{10}$ and $E_{13}$ together with Lemma~\ref{s2} imply that $m$ and $m'$
can play the role of opposite corners in $B$.

Assume that in the block structure of $G$, there are blocks/tracks $T_1,T_2,\dots T_l$ adjacent
to the same multiple node $m$.
Lemma~\ref{s3} and Lemma~\ref{s4} together with the excluded minors $E_5,E_5,E_7,E_8$
imply that $m$ is the only multiple node in $T_i$ and $m$ can be a degenerate side in $T_i$
with at most two exceptions.

These claims together prove the validity of $(iii)$.
\end{proof}

\section{Tree-like parts}

We sketch how to prove --- along the above ideas --- an excluded minor
characterization of PW2-graphs.
The result is not new, since the list has been obtained by Kinnersley and Langston \cite{Kinn94}
making a computer search.
Our proof outline explains certain similarities among the excluded minors.
Also, our proof will associate a task to each excluded minor (as it happened
in the 2-edge-connected case) that explains the role of that specific
excluded minor.

We know from basic graph theory that $G$ is built up from its blocks
and trees (consisting of edges, that are cut-edges of $G$) pasted together along
vertices in a tree-like fashion. 
A subgraph of $G$ which is a tree and consists only of
cut-edges (and hence it is an induced subgraph) is called {\sl a tree-part of $G$}.

One of the major reasons of difficulty in the discussion below 
is that these tree-like parts are not well-defined. 
A tree-like part can be considered as its tree-like edges are glued together. 
Although the maximal tree-like parts are uniquely defined,
we can not allow ourselves to consider the tree-like parts to be the maximal ones.

\begin{defn}
Let $T$ be a tree with two distinguished nodes $\ell$ and $r$.
We call $(T,\ell,r)$ a {\sl bond-tree} if and only if 
there exists a track $R\supseteq T$ such that the unique $\ell r$ path in $T$ can
be identified with a side of $R$.
The nodes $\ell$ and $r$ are opposite corners of the bond-tree.
\end{defn}

\begin{figure}[ht]
	\begin{center}
	\includegraphics[scale=0.5]{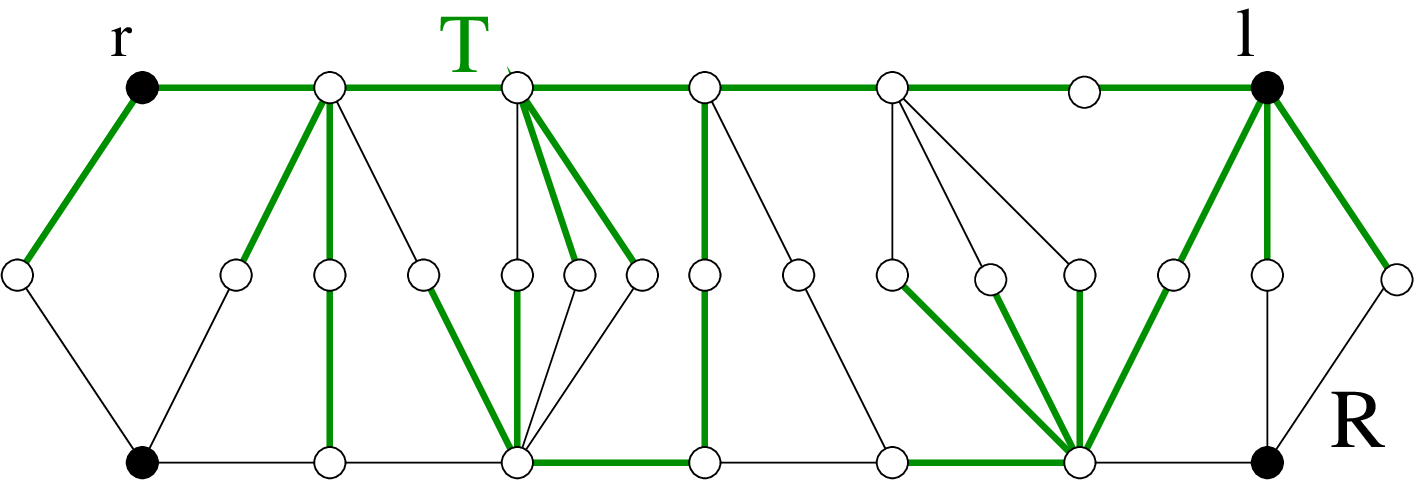}
 	\end{center}
  \caption{A bond-tree packed in a track}
\end{figure}

\begin{defn}
A subtree $T$ of $G$ is called a {\sl tree-frippery} rooted at $r$ if and only if 
\begin{enumerate}
\item[(a)]
the node $r$ is a cutnode of $G$, where
$T$ and $G\setminus T$ are glued together,
\item[(b)]
the subtree $T$ is a
subgraph of a track $R$ such that $T$ is disjoint from
one side of $R$, and $r$ can play the role of a corner.
\end{enumerate}
\end{defn}

A special kind of tree-frippery is important in the
structure theorem.

\begin{defn}
Let $(T,r)$ be a tree-frippery. 
It is called an {\sl edge-frippery} or {\sl hair}
if it has only one edge. 
Therefore, one of its endpoints has degree $1$,
and the other endpoint $r$ belongs to the rest of $G$.
\end{defn}

Blocks (tracks) and bond-trees can be glued together
as in the case of 2-edge-connected graphs.
In addition, there might be fripperies attached to the tracks.
The exact gluing pattern is described in the third
part of the following theorem.

\begin{thm}
Let $G$ be a connected graph.
The following three statements are equivalent.
\begin{enumerate}
\item[(i)]
The graph $G$ has path-width at most $2$.
\item[(ii)]
The graph $G$ does not contain any minor from the
Kinnersley--Langston list \cite{Kinn94}.
\item[(iii)]
The graph $G$ can be separated into its blocks and some tree parts, that
are classified as bond-trees, tree-fripperies and hairs.
The spine of $G$ can be glued together from
$2$-connected tracks (blocks) and bond-trees 
in a path-like fashion according to the following scheme.
The blocks and bond-trees are
enumerated as $B_1,B_2,\dots, B_L$ and there is a corresponding
sequence of nodes $m_1,m_2,\dots, m_L, m_{L+1}$ (consecutive nodes possibly coincide)
such that for any $2\le i\le L-1$ the nodes $m_i$ and $m_{i+1}$
are opposite corners of $B_i$.
Hence, $m_2, m_3,\dots, m_L$ are multiple nodes of the separation. 
The corners of blocks that are not multiple nodes are called {\sl free corners}.\\
There can be at most one tree-frippery attached to each free corner.\\
Finally, there might be hairs attached to the sides of the blocks.
So a hair can not be rooted at the middle node of a long
chord, but several hairs may be rooted at a side vertex of a track.
\end{enumerate}
\end{thm}

Before the proof sketch, we must clarify an important point.
Suppose there is a $B_i$ with degenerate side $\ell(=r)$ such that $\ell$ is a free corner.
That means $r$ is a free corner too, somewhat invisible.
There might be a tree-frippery attached to $\ell$, and another one attached to $r$.
Therefore, if a maximal tree part $T$ is attached to the track at
$\ell$, that is a free corner, then we must separate $T$ into a left and a right tree-frippery
in order to recognize the structure in (iii).
That is why the conditional form in (iii) is essential.

As before, the implications $(iii)\Rightarrow (i)\Rightarrow (ii)$ are straightforward.
The heart of the theorem is the implication $(ii)\Rightarrow (iii)$.
This can be proven by reversing the easy implications (see the original proof in \cite{Kinn94}
for the equivalence of (i) and (ii)) or by analyzing the cases when the gluing pattern 
in (iii) is not identifyable.
This proof method sheds light on the role of the excluded minors.
It is tedious, since the number of excluded minors is large.
The advantages of our method are already presented in the previous sections on
2-connected and 2-edge-connected graphs.
We do not want to fatigue the reader by the long case analysis, we only
want to exhibit the possibility of this alternative proof. 
This explains that the large number of excluded minors is
caused by the subtlety in the structure described in (iii),
due to the tree-like parts of $G$.

Let us demonstrate a few steps of the long road. 
The tree-like lemmas are without proofs.

\begin{lem} \label{hair}
Let $T$ be a tree with a distinguished node $r$.
The graph $T$ can not be a bucket of hairs rooted at the same node $r$
if and only if $T\not\succeq_r I$.
\end{lem}

\begin{figure}[ht]
	\begin{center}
	\includegraphics[scale=0.5]{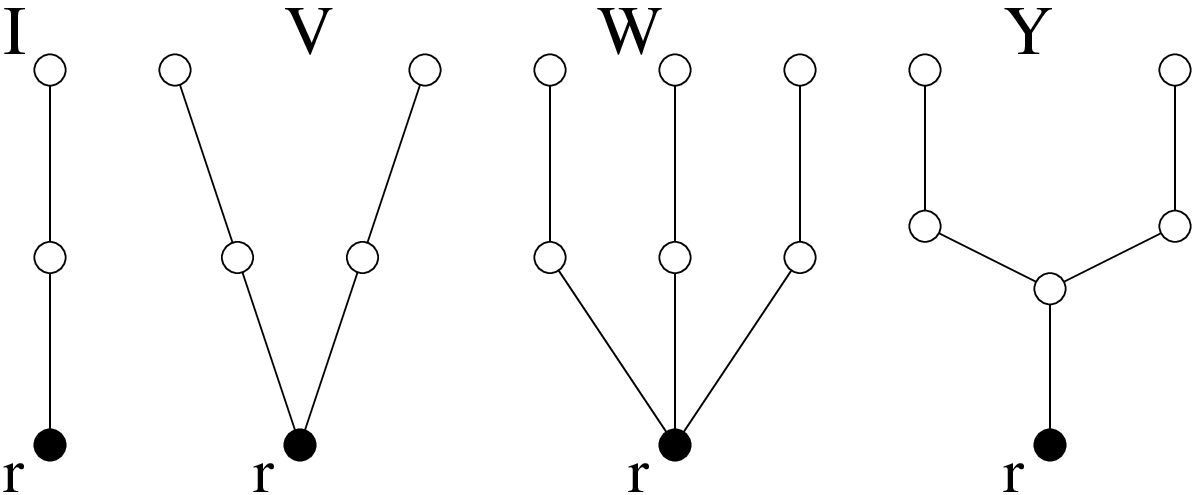}
 	\end{center}
  \caption{Attachments}
  \label{attach}
\end{figure}

\begin{lem} \label{fripp}
Let $T$ be a tree with a distinguished node $r$.
The graph $T$ can not be a tree-frippery rooted at $r$
if and only if $T\not\succeq_r V$. 
\end{lem}

\begin{lem} \label{fold}
Let $T$ be a tree with a distinguished node $r$.
The graph $T$ can not be a disjoint union of two tree-fripperies rooted at
the same node $r$ if and only if
$T\not\succeq_r W, Y$. 
\end{lem}

\begin{proof}
We outline a couple of steps.

\begin{itemize}
\item 
There must be a graph on the list of excluded minors guaranteeing that a
track can not have more than four points, where non-hair trees are
attached. The graph $I_1$ appears.
\item 
Assume that four trees are attached to a track.
There must be an excluded minor showing that at most two of them are not
tree-fripperies. See $I_4$.
\item 
Assume there are three nodes where the
attached trees are not tree-fripperies.
Then at least one of them must have an attachment that is a disjoint union 
of two tree-fripperies rooted at the same node.
This phenomenon is described
by the excluded minors $I_7, I_{10}, I_{12}, I_{14}$.
\item
We obtain a wider class of excluded minors if
we take into account the possible
connected tracks to our initial one (this fills the $I$-class in \cite{bar2}).
\item
Further classes arise when some inner structure is known about the
track where we want to see the neighborhood as described in (iii).
\item
We need more excluded minors to enforce the
right positions of hairs.
\item
Finally the pure tree case should be discussed.
\end{itemize}
\end{proof}

\begin{figure}[ht]
	\begin{center}
	\includegraphics[scale=0.5]{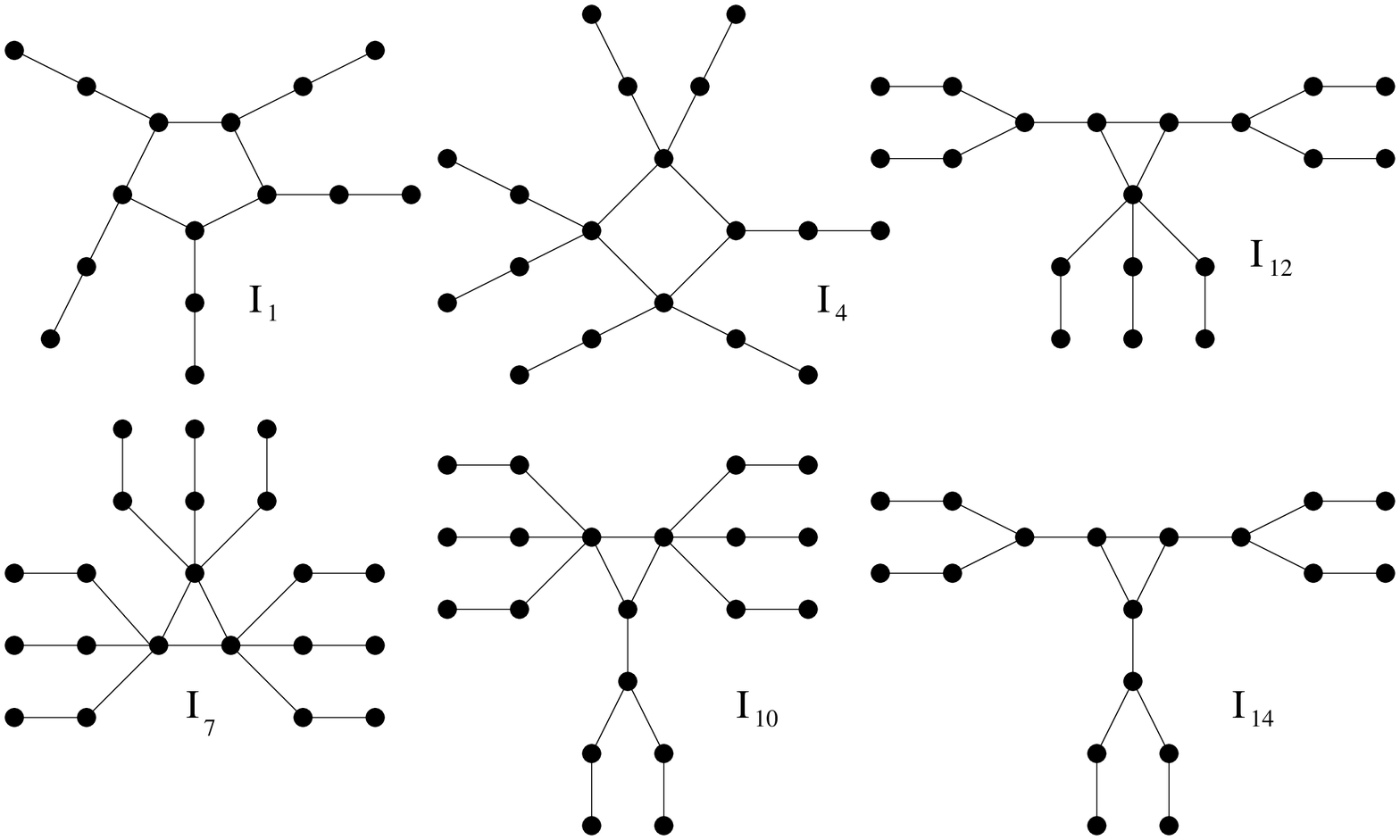}
 	\end{center}
  \caption{Illustrious examples of excluded minors}
  \label{napok}
\end{figure}


\section*{Conclusion}

We believe that $3$-connected PW3-graphs can be similarly characterized.
We achieved some results in this direction and plan to complete
those efforts. 

It looks substantially more difficult to continue along this line
to exact description of $k$-connected PW$k$-graphs.
We can expect asymptotic results rather than a precise one.
We raise the following questions:
\begin{itemize}
 \item What is the number of excluded minors for $k$-connected 
PW$k$-graphs in terms of $k$?
 \item Is there a lower bound, that is greater than polynomial in $k$?
 \item Is there a good upper bound?
\end{itemize}


\section*{Epilogue}

The first two authors initiated the study of PW2 in 1998,
when the first author started his PhD studies.
We were unaware of \cite{Kinn94} and looked for a characterization by paper and pencil.
The results we obtained were written in two article submissions in 1999.
Some version appeared in Bar\'at's Phd thesis \cite{phd}.
One article was published \cite{Barat01}, and the other was rejected \cite{bar2}.
While this other paper was rewritten, the last two authors submitted  their results \cite{lin}
for publication in 2003.
The first author was asked to referee the paper, and this is where our roads crossed.
The four authors agreed to unify their forces.
Due to various difficulties in space and time, the process lasted longer than it should have.
Finally, we agreed to publish the present version, which is seemingly rather different from
any of \cite{bar2} and \cite{lin}.
Therefore, we decided to leave those versions available on our home page.

\end{document}